\documentclass{birkmult}
\usepackage[a4paper,DIV=10,BCOR=20mm]{typearea}
 \newtheorem{theorem}{Theorem}[section]
 \theoremstyle{definition}
 \newtheorem{assumption}{Assumption}

 \theoremstyle{remark}

 \numberwithin{equation}{section}

 \usepackage{tikz,subfigure,amssymb,float}
 
\begin{document}
%
%
%
%
%
%
%
%
%
\title[Multiscale unique continuation properties of eigenfunctions]
 {Multiscale unique continuation properties of eigenfunctions}
\author[D.~Borisov]{Denis Borisov}
\address{ 
Institute of Mathematics, USC RAS \br
Chernyshevskii st., 112 \br
Ufa, 45008, \br
Russia
\br \& \br 
Bashkir State Pedagogical University \br 
October rev. st., 3a \br
Ufa, 450000 \br
Russia
}
\email{borisovdi@yandex.ru}

\author[I.~Naki\'c]{Ivica Naki\'c}
\address{
University of Zagreb \br 
Departement of Mathematics \br
10000 Zagreb \br
Croatia}
\email{nakic@math.hr}

\author[C.~Rose]{Christian Rose}
\address{
Chemnitz University of Technology \br
Faculty of Mathematics \br
09107 Chemnitz \br
Germany
}
\email{christian.rose@mathematik.tu-chemnitz.de}

\author[M.~Tautenhahn]{Martin Tautenhahn}
\address{
Chemnitz University of Technology \br
Faculty of Mathematics \br
09107 Chemnitz \br
Germany
}
\email{martin.tautenhahn@mathematik.tu-chemnitz.de}
%
%
\author[I.~Veseli\'c]{Ivan Veseli\'c}
\address{
Chemnitz University of Technology \br
Faculty of Mathematics \br
09107 Chemnitz \br
Germany
}
\email{ivan.veselic@mathematik.tu-chemnitz.de}
\thanks{D.B. was partially 
supported by Russian Science Foundation (project no. 14-11-00078) and the fellowship
of Dynasty foundation for young mathematicians.
I.N., Ch.R., M.T., and I.V. habe been partially supported by the DAAD and the Croatian Ministry of Science, Education and Sports
through the PPP-grant `Scale-uniform controllability of partial differential equations'.
M.T. and I.V. habe been partially supported by the DFG}


\date{\today}

\begin{abstract}
Quantitative unique continuation principles for multiscale structures are an important ingredient in a number applications, 
e.g.\ random Schr\"odinger operators and control theory. 

We review recent results and announce new ones regarding quantitative unique continuation principles for 
partial differential equations with an underlying multiscale structure. They concern Schr\"odinger and second order elliptic operators.
An important feature is that the estimates are scale free and with quantitative dependence on parameters. 
These unique continuation principles apply to functions satisfying certain `rigidity' conditions, 
namely that they are solutions of the corresponding elliptic equations, or projections on spectral subspaces. 
Carleman estimates play an important role in the proofs of these results. 
\end{abstract}

\maketitle
%
%
%
%
 \section{Introduction}
 \subsection*{Motivation: Retrival of global properties from local data}
In several branches of mathematics, as well as in applications, 
one often encounters problems of the following type:
Given a region in space $\Lambda \subset \mathbb{R}^d$, 
a subset $S \subset \Lambda $, 
and a function $f \colon \Lambda \to \mathbb{R}$,
what can be said about certain properties of $f \colon \Lambda \to \mathbb{R}$
given certain properties of $f|_S \colon S \to \mathbb{R}$?
In specific cases one may want to reconstruct $f$ as accurately as possible
based on knowledge of $f|_S$, in others it may be sufficient to estimate some
features of $f$.
\par
It is clear that for this task additional global information 
on $f$ is needed. Indeed, if $f$ is one of the indicator functions 
$\chi_S$ or  $\chi_{\Lambda\setminus S}$, an estimate based on $f|_S$
would yield wrong results. The first helpful property which comes to one's	 
mind is some regularity or smoothness property of  $f$. 
However, since  there are $C^\infty$-functions supported inside $S$ 
(or inside $\Lambda\setminus S$) this is not quite the right condition. 
The required property of $f$ is more adequately described as \emph{rigidity}, 
as we will see in specific theorems formulated below.

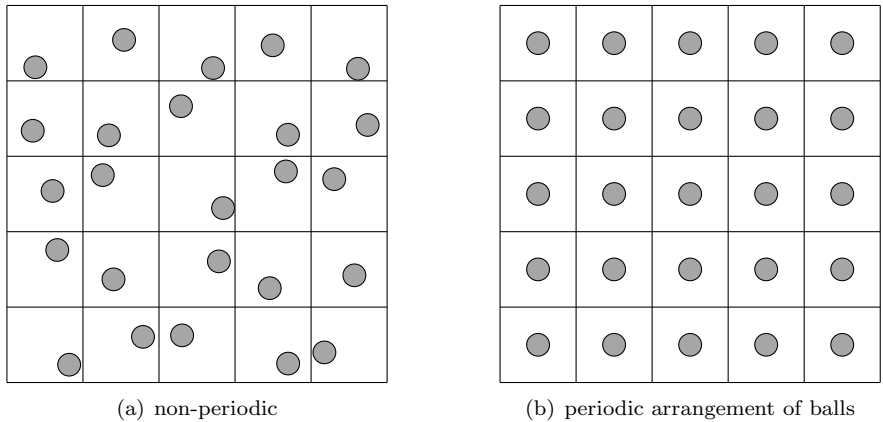
\begin{figure}[ht]\centering

\subfigure[non-periodic]{
\begin{tikzpicture}
\pgfmathsetseed{{\number\pdfrandomseed}}
\foreach \x in {0.5,1.5,...,4.5}{
  \foreach \y in {0.5,1.5,...,4.5}{
    \filldraw[fill=gray!70] (\x+rand*0.35,\y+rand*0.35) circle (0.15cm);
  }
}
\foreach \y in {0,1,2,3,4,5}{
\draw (\y,0) --(\y,5);
\draw (0,\y) --(5,\y);
}
\end{tikzpicture}
}
\hspace{1cm}
\subfigure[periodic arrangement of balls]{
\begin{tikzpicture}
\foreach \x in {0.5,1.5,...,4.5}{
  \foreach \y in {0.5,1.5,...,4.5}{
    \filldraw[fill=gray!70] (\x,\y) circle (1.5mm);
  }
}
\foreach \y in {0,1,2,3,4,5}{
  \draw (\y,0) --(\y,5);
  \draw (0,\y) --(5,\y);
}
\end{tikzpicture}
}
\caption{Examples of equidistributed sets $S$ within region $\Lambda\subset \mathbb{R}^2$.\label{fig:equidistributed}}
\end{figure}

In this paper we are mainly concerned with problems with a multiscale structure.
For this reason it is natural to require that the set $S$ is in some sense equidistributed 
within $\Lambda$.  
At this point we will not give a precise definition of such sets.
It will become clear that such a set $S$ should be relatively dense in $\mathbb{R}^d$ or $\Lambda$,
and should have positive density. A particularly nice set $S$ would be a 
periodic arrangement of balls, and we want to include small perturbations of such a configuration.
Thus, equidistributed sets could be seen as a generalization of such a situation, cf.\ Fig.~\ref{fig:equidistributed}.

\subsection*{Example: Shannon sampling theorem}
We recall a well known theorem as an example
or benchmark, see e.g.\ \cite{ButzerSS-88}. This way we will see what is the best we can hope for   
in the task of reconstructing a function. Moreover, we will encounter 
one possible interpretation what the term \emph{rigidity} means, 
and see major differences between the reconstruction problem in dimension one and higher dimensions.
\par
The Shannon sampling theorem states: 
Let $f\in C(\mathbb{R}) \cap L^2(\mathbb{R})$ be such that the Fourier transform 
\[
 \hat f (p)= \frac{1}{\sqrt{2\pi}} \int_\mathbb{R} e^{-i \, x \, p }\, f(x) \, \mathrm{d} x 
\]
vanishes outside $[-\pi K, \pi K]$.
Then the series 
\begin{equation}
\label{eq:Shannon}
 (S_K f)(x) = \sum_{j \in \mathbb{Z}} f\left( \frac{j}{K}\right) \frac{\sin\pi (K \, x - j)}{\pi (K \, x - j)}
\end{equation}
converges absolutely and uniformly for $x \in \mathbb{R}$ and 
\[
 S_K f=f\ \text{on}\ \mathbb{R} .
\]
Thus we can reconstruct the original function $f$ from the sample values  $ f( j / K)$,
which are multiplied with weights depending on the distance to the point $x \in \mathbb{R}$ and summed up.
Here the rigidity condition is implemented by the requirement $\operatorname{supp} \hat f \subset [-\pi K, \pi K]$, which implies that $f$ is entire.
A remarkable feature of this exact result is that it is stable under perturbations:
If the nodes $j$ deviate slightly from the integers, or if the measurement data 
$ f( \frac{j}{K})$ are inaccurate, the error $f-S_K f$ can still be controlled.
If the support condition $\operatorname{supp} \hat f \subset [-\pi K, \pi K]$ is violated, the \emph{aliasing error} is estimated as
\begin{equation}
\label{eq:aliasing}
 \sup_{x\in\mathbb{R}} \lvert f(x)-S_K f(x) \rvert \leq \sqrt{\frac{2}{\pi}} \int_{\lvert p \rvert >\pi K} \lvert \hat f(p) \rvert \, \mathrm{d} p .
\end{equation}
This will give, for instance, good results for centered Gaussians with appropriate variance.

Statements \eqref{eq:Shannon} and \eqref{eq:aliasing} are  strong with respect to the sampling set $S=\mathbb{Z}$, which is very thin. It has zero Lebesgue measure, in fact, it is discrete.
Albeit, it is relatively dense in $\mathbb{R}$, so it has some of the properties we associated with an equidistributed set.
Compared to Shannon's theorem, the results we present below appear  much weaker.
This is, among others, due to two features: we consider functions on multidimensional space,
which, in addition, have low regularity, in fact are defined as equivalence classes in some $L^2$
or Sobolev space. In this situation evaluation of a function at a point may not have a proper meaning.
This is one of the reasons why we have to consider samples $S$ which are composed of small balls, rather than  single points. 
A second aspect where dimensionality comes into play is the following: 
A polynomial of one variable of degree $N$ vanishes identically if it has $N+1$ zeros.  
A non-trivial polynomial in two variables may vanish on an uncountable set (albeit not on one of positive measure).
This illustrates that reconstruction estimates for functions of several variables are more subtle than Shannon's theorem. 
Consequently, one has to settle for more modest goals than the full reconstruction of the function $f$. We want to  derive an equidistribution property for functions satisfying some 
rigidity property. As will be detailed later this result is called --- depending on the context and scientific environment --- scale free unique continuation property, observability estimate, or uncertainty relation.
A first result of this type is formulated in the next section.
\section{Equidistribution property of Schr\"odinger eigenfunctions}
The following result \cite{Rojas-MolinaV-13} was motivated by questions arising in the spectral theory of random Schr\"odinger operators. Later, it turned out that similar estimates are of relevance in the control theory of the heat equation.

We fix some notation. For $L>0$ we denote by $\Lambda_L = (-L/2 , L/2)^d$ a cube in $\mathbb{R}^d$.
For $\delta>0$ the open ball centered at $x\in \mathbb{R}$ with radius $\delta$ is denoted by 
$B(x, \delta)$. For a sequence of points $(x_j)_j$ indexed by $j \in \mathbb{Z}^d$
we denote the collection of balls  $\cup_{j \in \mathbb{Z}^d} B(x_j , \delta) $
by $S$ and its intersection with $\Lambda_L$ by $S_L$. 
%
%
We will be dealing with certain self-adjoint operators on subsets of $\mathbb{R}^d$.
Let  $\Delta$ be the $d$-dimensional Laplacian, $V \colon \mathbb{R}^d\to \mathbb{R}$ a bounded measurable
function, and $H_L = (-\Delta + V)_{\Lambda_L} $ a Schr\"odinger operator on the cube 
$\Lambda_L$ with Dirichlet or periodic boundary conditions.
The corresponding domains are denoted by 
$\mathcal{C}(\Delta_{\Lambda,0})\subset W^{2,2}(\Lambda_L)$ and $ \mathcal{C}(\Delta_{\Lambda,\mathrm{per}})$, respectively. 
Note that we denote a multiplication operator by the same symbol as the corresponding function. 

\begin{theorem}[\cite{Rojas-MolinaV-13}]  
\label{thm:RojasVeselic}
 Let $\delta, K > 0$. Then there exists $C \in (0,\infty)$ such that for all
 $L \in 2\mathbb{N}+1 $, all measurable   $V : \mathbb{R}^d \to [-K , K]$, all real-valued $\psi \in 
 \mathcal{C}(\Delta_{\Lambda,0})  \cup  \mathcal{C}(\Delta_{\Lambda,\mathrm{per}})$
 with $(-\Delta + V)\psi = 0$ almost everywhere on $\Lambda_L$, and  all sequences $(x_j)_{j \in \mathbb{Z}^d} \subset \mathbb{R}^d$, 
 such that  $\forall j \in \mathbb{Z}^d: \, B(x_j , \delta) \subset \Lambda_1 + j$ we have
\begin{equation}
\label{eq:observability}
\int_{S_L} \psi^2  \geq C \int_{\Lambda_L} \psi^2 .
\end{equation}
\end{theorem} 

To appreciate the result properly, the quantitative dependence of the constant $C$ on model parameters
is crucial. The very formulation of the  theorem  states that $C$ is
independent of position of the balls $B(x_j,\delta)$ within $\Lambda_1 +j$, 
and independent of the scale $L\in2\mathbb{N} +1$.
The estimates given in Section 2 of \cite{Rojas-MolinaV-13} show moreover, 
that $C$ depends on the potential
$V$ only through the norm $\lVert V \rVert_\infty$ (on an exponential scale), 
and it depends on the small radius $\delta>0$ polynomially, i.e.\ $C\gtrsim \delta^N$, for some $N\in\mathbb{N}$ which depends on the dimension $d$ and $\lVert V \rVert_\infty$. 
This shows that we are not able to control the integral $\int_{S_L} \psi^2 $
by evaluating $\psi$ at the midpoints $j \in \mathbb{Z}^d$ of the unit cubes.
One sees with what rate the estimate diverges, as the balls become smaller and approximate a single point.
The polynomial behavior $C\gtrsim \delta^N$ can be readily understood when looking at monomials $\psi_n(x)=x^n$ on the unit interval $(0,1)$. There we have
\[
\int_{(0,\delta)} \psi_n^2= \frac{\delta^{2n+1}}{2n+1}
=\delta^{2n+1} \int_{(0,1)} \psi_n^2 .
\]
We formulated the theorem only for the eigenvalue zero, but it is easily applied to other eigenfunctions as well since
\[
 H_L\psi=E\psi \Leftrightarrow (H_L-E)\psi=0 .
\]
Consequently the constant $K=K_{V}$ has to be replaced 
with the possibly larger $K=K_{V-E}$.

There is a very natural question, which was spelled out in \cite{Rojas-MolinaV-13}, 
namely does the same estimate  \eqref{eq:observability} hold 
true for linear combinations $\psi\in \operatorname{Ran} \chi_{(-\infty,E]} (H_L)$ of eigenfunctions as well?
The property in question can be equivalently stated as:
Given $\delta >0, K\geq0, E\in\mathbb{R}$ there is a constant $C>0$ such that 
for all measurable $ V\colon \mathbb{R}^d \rightarrow [-K,K] $, all $L \in 2\mathbb{N}+1$,
and all sequences $(x_j)_{j\in\mathbb{Z}^d} \subset \mathbb{R}^d$ 
with $B(x_j,\delta)  \subset\Lambda_1 +j$ for all $j \in \mathbb{Z}^d$ we have
\begin{equation}
\label{eq:uncertainty}
 \chi_{(-\infty,E]} (H_L) \, W_L \, \chi_{(-\infty,E]} (H_L) \geq C~ \chi_{(-\infty,E]} (H_L) ,
\end{equation}
where $W_L=\chi_{S_L}$  is the indicator function of  $S_L$ and 
$\chi_{I} (H_L)$ denotes the spectral projector of $H_L$ onto the interval $I$.
Here $C=C_{\delta, K, E}$ is determined by $\delta, K, E$ alone.

Note that all considered operators are lower bounded by $-K$ in the sense of quadratic forms.
Thus the spectral projection on the energy interval $(-\infty,E]$ is the same as the
spectral projection on the energy interval $[-K,E]$. The upper bound $E$ in the energy parameter
is crucial for preventing the corresponding eigenfunctions to oscillate too much.
 
One can pose a modified version of the question: 
Given $\delta >0, K\geq0, a<b\in\mathbb{R}$ is there is a constant $\tilde C>0$ 
such that for all measurable $ V\colon \mathbb{R}^d \rightarrow [-K,K] $, all $L \in 2\mathbb{N}+1$,
and all sequences $(x_j)_{j\in\mathbb{Z}^d} \subset \mathbb{R}^d$ 
with $B(x_j,\delta)  \subset\Lambda_1 +j$ for all $j \in \mathbb{Z}^d$ we have
\begin{equation}
\label{eq:local_uncertainty}
 \chi_{[a,b]} (H_L) \, W_L \, \chi_{[a,b]} (H_L) \geq \tilde C~ \chi_{[a,b]} (H_L) .
\end{equation}
Here $\tilde C=\tilde C_{\delta, K, a, b}$ depends (only) on $\delta, K, a, b$.
Note that  inequality \eqref{eq:uncertainty} implies
\eqref{eq:local_uncertainty} since
\begin{align*}
\chi_{[a,b]} (H_L) \, & W_L \, \chi_{[a,b]} (H_L) \\[1ex]
&=\chi_{[a,b]} (H_L) \, \chi_{(-\infty,b]} (H_L) \, W_L \, \chi_{(-\infty,b]} (H_L) \, \chi_{[a,b]} (H_L) 
\\[1ex]
&\geq C_{\delta, K, b}~ \chi_{[a,b]} (H_L) \, \chi_{(-\infty,b]} (H_L) \, \chi_{[a,b]} (H_L) 
\\[1ex] &=C_{\delta, K, b}  \chi_{[a,b]} (H_L) .
\end{align*}
However, $C_{\delta, K, b}$ may be substantially  smaller than
$\tilde C_{\delta, K, a, b}$
due to the enlarged energy interval.

Klein obtained a positive answer to the question for sufficiently short intervals. 
\begin{theorem}[\cite{Klein-13}] \label{thm:Klein-13}
Let $d \in \mathbb{N}$, $E\in \mathbb{R}$, $\delta\in (0,1/2]$ and $V:\mathbb{R}^d \to \mathbb{R}$ be measurable and bounded. There is a constant $M_d>0$ such that if we set
\[
 \gamma = \frac{1}{2} \delta^{M_d \bigl(1 + (2\lVert V \rVert_\infty + E)^{2/3}\bigr)} ,
\]
then for all energy intervals $I\subset (-\infty, E]$ with length bounded by $2\gamma$, all $L \in 2\mathbb{N}+1$, $L\geq 72 \sqrt{d}$
and all sequences $(x_j)_{j\in\mathbb{Z}^d} \subset \mathbb{R}^d$ 
with $B(x_j,\delta)  \subset\Lambda_1 +j$ for all $j \in \mathbb{Z}^d$
\begin{equation}
 \chi_{I} (H_L) \, W_L \, \chi_{I} (H_L) \geq \gamma^2\chi_{I} (H_L) .
\end{equation}
\end{theorem}
Although this does not answer the above posed question 
for arbitrary compact intervals, the result is sufficient for many questions in spectral theory of random Schr\"o\-dinger operators. 
A generalization of Theorem~\ref{thm:Klein-13} to intervals of arbitrary length is given in Section~\ref{sec:NTV}.
This answers completely the question posed in \cite{Rojas-MolinaV-13}.
\par
Depending on the context and the area of mathematics 
the above described estimates carry various names. 
If one speaks of an \emph{equidistribution property of eigenfunctions}, one is interested in the 
comparison of the measure $\lvert \psi(x) \rvert^2 \mathrm{d}x$ with the uniform distribution
on the cube $\Lambda_L$.
The term \emph{scale free unique continuation principle} is used in works concerning random Schr\"odinger 
operators. It refers to a quantitative version of the classical unique continuation principle, 
which is uniform on all large length scales.
One can interpret Theorem~\ref{thm:RojasVeselic} as an \emph{uncertainty relation}:
the condition $H_L\psi=E \psi $ corresponds to a restriction in momentum/Fourier-space 
and enforces a delocalization/flatness property in direct space.
Similarly, the spectral projector 
$\chi_{(-\infty,E]}$ in Ineq.~\eqref{eq:uncertainty} corresponds to a restriction in momentum space.
Here we see a direct analogy to Shannon's theorem discussed above: If the Fourier transform of a function is sufficiently concentrated, the function itself cannot vary too much over short distances.
Ineq.~\eqref{eq:local_uncertainty} can also be interpreted as a \emph{gain of positive definiteness}. It says that for a general self-adjoint operator $A\geq 0$, which may have a kernel, 
and an appropriately chosen spectral projector $P$ of the Hamiltonian, the restriction  $P A P \geq c P$ is strictly positive.
In control theory results as we discuss them are sometimes called \emph{observability} estimates.
This term is more common for time-dependent partial differential equations, but sometimes used for stationary ones as well.
\par
In the literature on random Schr\"odinger operators related results have been
derived before in a number of papers. For more details we refer to Section 1 of \cite{Rojas-MolinaV-13}.


\section{Methods and background}
A paradigmatic result for the \emph{weak unique continuation principle} is the following. A solution of $ \Delta f \equiv 0 \text{~on~} \mathbb{R}^d $ satisfying  
$f\equiv0\ \text{on}\ B(0,\delta)$ for arbitrary small, but positive $\delta$, must vanish on all of
$\mathbb{R}^d$. 
The restrictive conditions can be relaxed. First of all, the condition  
$f\equiv0 \text{~on~} B(0,\delta)$ can be replaced by 
\begin{displaymath}
 \forall N\in\mathbb{N} \qquad\lim\limits_{\delta\searrow0} \delta^{-N} \int\limits_{B(0,\delta)} \lvert f (x) \rvert \mathrm{d} x=0 .
\end{displaymath}
In this form the implication is called \emph{strong unique continuation principle}.
Moreover, the Laplacian $\Delta $ can be replaced by a rather general second order elliptic operator.
We will discuss related results in Sections \ref{sec:NRT} and \ref{sec:BTV}. 
A powerful method to prove unique continuation statements, as well as quantitative versions thereof, 
are Carleman estimates. Originally, Carleman \cite{Carleman-39} derived them for functions of two variables. 
Later M\"uller \cite{Mueller-54} extended the estimates to higher dimensions.
By now, there are hundreds of papers dealing with Carleman estimates.
We will describe one explicit version in Section~\ref{sec:NRT}, which is an important tool
for the quantitative unique continuation estimates discussed shortly for Schr\"odinger operators.
In Section~\ref{sec:BTV} we will present new results in this direction which deal with elliptic 
second order operators with variable coefficients.

\subsection*{Quantitative unique continuation principle}

In \cite{BourgainK-05} Bourgain and Kenig derived the following pointwise 
quantitative unique continuation principle.
\begin{theorem}
 Assume $(-\Delta +V) u= 0$ on $\mathbb{R}^d$ and $u(0) = 1$, 
 $\lVert u \rVert_\infty \leq C$, $\lVert V \rVert_\infty \leq C$. Let $x_0 \in \mathbb{R}^d$, $\lvert x_0 \rvert = R > 1$. 
 Then there exists a constant $C' > 0$ such that
\[
\max_{\lvert x-x_0 \rvert \leq 1} \lvert u(x) \rvert > C' \exp \left( -C' (\log R)R^{4/3} \right)  .
\]
\end{theorem}


In our context a version of this result with local $L^2$-averages is more appropriate.
Various estimates of this type have been given in \cite{GerminetK-13,BourgainK-13,Rojas-MolinaV-13}. 
We quote here the version from the last mentioned paper.

\begin{theorem}
\label{thm:quantitative-UCP}
Let $K, R, \beta\in [0, \infty), \delta \in (0,1]$.
There exists a constant $C_{\rm qUC}=C_{\rm qUC}(d,\allowbreak{}K_V, R,\delta, \beta) >0$ such that,
for any $G\subset  \mathbb{R}^d$ open, any $\Theta\subset G$ measurable,
satisfying the geometric conditions
\[
\operatorname{diam} \Theta +  \operatorname{dist} (0 , \Theta) \leq 2R \leq 2  \operatorname{dist} (0 , \Theta), \quad \delta < 4R, \quad B(0, 14R ) \subset G,
\]
and any measurable $V\colon G \to [-K,K]$  and real-valued $\psi\in W^{2,2}(G)$ satisfying the  differential inequality
\begin{equation*}
\label{eq:subsolution}
\lvert \Delta \psi \rvert \leq \lvert V\psi \rvert \quad \text{a.e. on }   G
\quad  \text{ as well as } \quad 
\int_{G}  \lvert \psi \rvert^2
\leq
\beta \int_{\Theta} \lvert \psi \rvert^2 ,
\end{equation*}
 we have
\begin{equation*}
\label{eq:aim}
\int_{B(0,\delta)}  \lvert \psi \rvert^2
\geq
C_{\rm qUC}
\int_{\Theta} \lvert \psi\rvert^2 .
\end{equation*}
\end{theorem}
\begin{figure}[ht]\centering
 \begin{tikzpicture}
  \filldraw[black!10] plot[smooth cycle] coordinates{(1,1) (1,6) (4,6) (8,7.5) (10,7) (9,2) };
  \draw plot[smooth cycle] coordinates{(1,1) (1,6) (4,6) (8,7.5) (10,7) (9,2) };
  \filldraw[black!30] (5.5,3.5) circle (0.5cm);
  \draw  (5.5,3.5) circle (0.5cm);
  \filldraw (5.5,3.5) circle (1pt);
  \draw (4.4,3.7) node {\small $B(0,\delta)$};
  \filldraw[black!30]   (7,4) rectangle (8,5);
   \draw  (7,4) rectangle (8,5);
  \draw (8.3,4.5) node {\small $\Theta$};
  \draw[latex-latex] (5.51,3.51)--(7,4.1);
  \draw[latex-latex] (5.5,3.48)--(5.8,1.095);
  \draw (6.35,4.1) node {\small $R$};
  \draw (5.3,2) node {\small $14R$};
  \draw (4,6.3) node {\small $G$};
\end{tikzpicture}
\caption{Assumptions in Theorem \ref{thm:quantitative-UCP} on the geometric constellation of $G$, $\Theta$, and $B(0,\delta)$}
\end{figure}
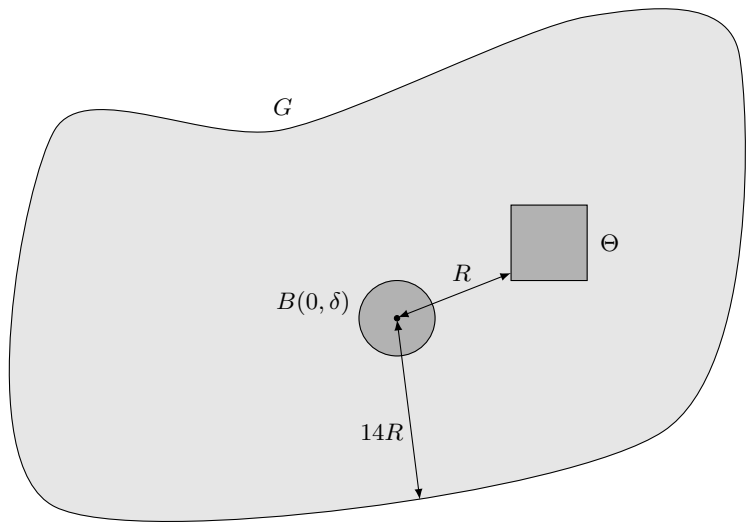


\section{Equidistribution property of linear combinations of eigenfunctions}\label{sec:NTV}

In this section we present a result from a project of I.~Naki\'c, M.~T\"aufer, M.~Tautenhahn and I.~Veseli\'c \cite{NakicTV}, namely which gives Ineq.~\eqref{eq:observability} also for linear combinations of eigenfunctions $\psi \in \operatorname{Ran} \chi_{(-\infty , E]} (H_L)$ for arbitrary $E \in \mathbb{R}$. As shown above, this implies Ineq.~\eqref{eq:uncertainty} for arbitrary $E \in \mathbb{R}$ and hence Ineq.~\eqref{eq:local_uncertainty} for $[a,b] \subset (-\infty , E]$. Indeed, our result gives a full answer to the open question in \cite{Rojas-MolinaV-13} whether Theorem~\ref{thm:quantitative-UCP} holds also for linear combinations of eigenfunctions, which was partially answered in \cite{Klein-13}, cf.\ Theorem~\ref{thm:Klein-13}. 

Since we first show Ineq.~\eqref{eq:uncertainty} for arbitrary $E\in\mathbb{R}$, the constant $\tilde C$ in Ineq.~\eqref{eq:local_uncertainty} will not be optimal, since it does not depend on the lower bound $a$ of the interval $[a,b] \subset [-\infty , E)$.

The following theorem was given in \cite{NakicTV} and full proofs will be
provided in \cite{NakicTTV-prep}.
\begin{theorem}[\cite{NakicTV}] \label{thm:NTV2013}
There is $N = N(d)$ such that for all $\delta \in (0,1/2)$, all measurable and
bounded $V: \mathbb{R}^d \to \mathbb{R}$, all $L \in \mathbb{N}$, all $E \geq 0$
and all $\psi \in \operatorname{Ran} (\chi_{(-\infty,E]}(H_L))$ and all
sequences $(x_j)_{j\in\mathbb{Z}^d}\subset \mathbb{R}^d$, such that for all
$j\in \mathbb{Z}^d$ $B(x_j,\delta)\subset \Lambda_1+j$, we have
\begin{equation} \label{eq:result1}
\int_{S_L} \lvert \psi \rvert^2 \geq C_{\mathrm{sfuc}}
\int_{\Lambda_L}\lvert \psi \rvert^2
\end{equation}
\[
C_{\mathrm{sfuc}} = C_{\mathrm{sfuc}} (d, \delta ,E,\lVert V \rVert_\infty )
:= \delta^{N \bigl(1 + \lVert V \rVert_\infty^{2/3} + \sqrt{E} \bigr)}.
\]
\end{theorem}
Hence, as in Theorem~\ref{thm:RojasVeselic}, the constant is independent on
the position of the balls $B(x_j,\delta)$, the scale $L$, and it depends on
the potential $V$ only through the norm $\lVert V\rVert_{\infty}$.

Here we give a sketch of the proof. 
We use two different Carleman inequalities in $\mathbb{R}^{d+1}$, one with a boundary term in $\mathbb{R}^d \times \{ 0 \}$ and the other without boundary terms. From these Carleman estimates we deduce two interpolation inequalities for a solution of a Schr\"odinger equation in $\mathbb{R}^{d+1}$. In the final step we apply these interpolation inequalities to the function $F : \Lambda_L \times \mathbb{R}
\to \mathbb{C}$ defined by
\begin{equation*} \label{eq:F}
F (x) = \sum_{\genfrac{}{}{0pt}{2}{k \in \mathbb{N}}{E_k \leq E}} \alpha_k
\psi_k (x’) s_k( x_{d+1}) ,
\end{equation*}
where $\alpha_k = \langle \psi_k , \psi \rangle$ with $\psi_k$ denoting the
eigenfunctions of $H_L$ corresponding to the eigenvalues $E_k$,
$\mathbb{R}^{d+1}\ni x=(x’,x_{d+1})$, $x’\in\mathbb{R}^d$, $x_{d+1}\in
\mathbb{R}$ and
\[
s_k(x)=\begin{cases}
     \sinh(\sqrt{E_k} x)/ \sqrt{E_k}, & E_k>0,\\
     x, & E_k=0,\\
     \sin(\sqrt{\lvert E_k \rvert} x)/ \sqrt{\lvert E_k \rvert}, & E_k<0 .
\end{cases}
\]
This function $F$ satisfies $\Delta F  = V F$ on $\Lambda_L \times \mathbb{R}$ and $\partial_{d+1} F (x',0) = \psi (x')$ on $\Lambda_L$, 
and one can obtain upper and lower estimates for the $H^1$-norm of the function $F$ in terms of the parameters $K$, $E$, $d$ and $\sum_{E_k \leq E}\lvert \alpha_k \rvert^2$.
\section{Explicit Carleman estimates for elliptic operators}\label{sec:NRT}
As mentioned above, Carleman estimates play a significant role in the results about unique continuation principles. In the case of quantitative unique continuation principles on multiscale structures, it is important to have a Carleman estimate with dependence on various parameters as precise as possible.

We consider the second order elliptic partial differential operator 
\[
L =-\sum_{i,j=1}^d \partial_i\left(a^{ij}\partial_j\right),
\]
acting on functions in $\mathbb{R}^d$. We introduce the following assumption on the coefficient functions $a^{ij}$.

\begin{assumption}\label{ass:elliptic+}
Let $r,\vartheta_1 , \vartheta_2 > 0$. The operator $L$ satisfies $A(r,\vartheta_1 , \vartheta_2)$, if and only if $a^{ij} = a^{ji}$ for all $i,j \in \{1,\ldots , d\}$ and for almost all $x,y \in B(r)$ and all $\xi \in \mathbb{R}^d$ we have
\begin{equation*} \label{eq:elliptic}
\vartheta_1^{-1} \lvert \xi \rvert^2 \leq \sum_{i,j=1}^d a^{ij} (x) \xi_i \xi_j \leq \vartheta_1 \lvert \xi \rvert^2 \quad\text{and}\quad \sum_{i,j=1}^d \lvert a^{ij} (x) - a^{ij} (y) \rvert \leq \vartheta_2 \lvert x-y \rvert .
\end{equation*}
\end{assumption}

Here $B (r) \subset \mathbb{R}^d$ denotes the open ball in $\mathbb{R}^d$ with radius $r$ and center zero.
Let the entries of the inverse of the matrix $( a^{ij}(x) )_{i,j=1}^d$ be denoted by $a_{ij}(x)$.

We present the result for the ball $B(1)$, but by scaling arguments this result can be generalized to arbitrary large balls $B(R)$, now with a different weight function which depends also on $R$.
\par
In the following theorem we formulate a Carleman estimate for elliptic
partial differential operators with variable coefficients analogous to
those given in \cite{EscauriazaV-03} for parabolic operators. In the
case of the pure Laplacian this has already been done in
\cite{BourgainK-05}. In particular, we establish that the estimate is
valid on the whole domain (i.e.\ $\delta =1$ holds in the notation of
\cite{EscauriazaV-03}) and give quantitative estimates for all the
parameters. This is part of a recent work of I.~Naki\'c, C.~Rose and
M.~Tautenhahn \cite{NakicRT}.
\par
For $\mu >0$ let $\sigma : \mathbb{R}^d \to [0,\infty)$ and $\psi
\colon [0,\infty)\to[0,\infty)$ be given by
\[
\sigma(x) = \left( \sum_{i,j=1}^{d} a_{ij}(0) x_i x_j \right)^{1/2}
\quad \text{and} \quad
\psi(s) = s \cdot \exp\left[-\int_0^s\frac{1- \mathrm{e}^{-\mu t}}t
\mathrm{d}t\right] .
\]
We define the weight function $w\colon\mathbb{R}^d\to [0,\infty)$ by
$w(x) = \psi(\sigma(x))$. Note that the weight function satisfies the
bounds
\begin{equation} \label{eq:bounds_w}
\forall x \in B(1) \colon \quad \frac{\lvert x \rvert}
{C_3 \sqrt{\vartheta_1}} \leq w(x) \leq \sqrt{\vartheta_1}
\lvert x \rvert \quad \text{with} \quad C_3 = \mathrm{e} \mu .
\end{equation}
\begin{theorem}[\cite{EscauriazaV-03,NakicRT}]
\label{thm:Carleman}
Let $\vartheta_1 , \vartheta_2 > 0$ and Assumption $A(1,\vartheta_1 ,
\vartheta_2)$ be satisfied. Then there exist constants $\mu,C_1,
C_2>0$ depending only on $\vartheta_1,\vartheta_2$ and the
dimension $d$ such that for all $f\in C_0^\infty(B(0,1)
\setminus\{0\})$ and all $\alpha>C_1$ we have
\[
\int \alpha w^{1-2\alpha}\lvert\nabla f\rvert^2+\alpha^3
w^{-1-2\alpha}f^2\leq C_2\int w^{2-2\alpha}(Lf)^2 \text.
\]
\end{theorem}
Explicit bounds on $\mu=\mu(\vartheta_1,\vartheta_2)$ are given in
\cite{NakicRT}. In particular,
\begin{equation}
 \label{eq:uniform-mu}
\forall \, T >0 : \quad \mu_T = \sup\{\mu(\vartheta_1,\vartheta_2)\colon 0<\vartheta_1,\vartheta_2\leq T\} < \infty .
\end{equation}
With a regularization procedure (see, for example, \cite[Theorem 1.6.1]{Ziemer-89}) this result can be extended to the functions in $H_0^2(B(0,1))$ which are compactly supported away from the origin. 

\section{Quantitative unique continuation estimates for elliptic operators}\label{sec:BTV}

In this section we announce a result from an ongoing work of D.~I.~Boris, M.~Tautenhahn and I.~Veseli\'c \cite{BorisovTV}. It concerns a quantitative unique continuation principle for elliptic second order partial differential operators with slowly varying coefficients. 

As in the previous section we denote by $L$ the second order partial differential operator 
\[
 Lu = -\sum_{i,j=1}^d \partial_i \left( a^{ij} \partial_j u \right) ,
\]
acting on functions $u$ on $\mathbb{R}^d$.

\begin{theorem}[\cite{BorisovTV}] \label{thm:qUC_elliptic} 
Let $R,\vartheta_1 , \vartheta_2 \in (0,\infty)$, $D_0 < 6 R$, $K_V, \beta \in [0,\infty)$, $\delta \in (0, 4 R]$, let $C_3 = C_3 (d,\vartheta_1,\vartheta_2)$ be the constant from Eq.~\eqref{eq:bounds_w}, and assume that 
\begin{equation*} \label{eq:slowly}
 A(12R+2D_0, \vartheta_1 , \vartheta_2 ) \quad \text{and} \quad \vartheta_1 C_3 < \frac{1}{4R}
\end{equation*}
are satisfied. Then there exists $C_{\rm qUC} = C_{\rm qUC} (d , \vartheta_1 , \vartheta_2 , R , D_0 , K_V , \delta , \beta) > 0$, such that, for any $G\subset \mathbb{R}^d$ open, $x \in G$ and $\Theta \subset G$ measurable, satisfying
 \[
 \operatorname{diam} \Theta + \operatorname{dist} (x , \Theta) \leq 2R \leq 2 \operatorname{dist} (x , \Theta) \quad
 \text{and} \quad
 B(x,12R+2D_0) \subset G,
\]
and any measurable $V : G \to [-K_V , K_V]$ and real-valued $\psi \in W^{2,2} (G)$ satisfying the differential inequality
\begin{equation*} \label{eq:psi}
 \lvert L \psi \rvert \leq \lvert V\psi \rvert \quad \text{a.e.\ on $G$} \quad \text{as well as} \quad \int_G \lvert \psi \rvert^2 \leq \beta \int_\Theta \lvert \psi \rvert^2 ,
\end{equation*}
we have
\[
 \int_{B(x,\delta)} \lvert \psi \rvert^2 \geq C_{\rm qUC} \int_{\Theta} \lvert \psi \rvert^2 .
\]
\end{theorem} 

Theorem~\ref{thm:qUC_elliptic} generalizes Theorem~\ref{thm:RojasVeselic} to second order elliptic operators with slowly varying coefficient functions. 
This is explicitly given by the assumption $\vartheta_1 C_3 < 1/ (4R)$. 
Indeed, for fixed $R>0$ the last inequality is satisfied for $\vartheta_1$ sufficiently small, since
\eqref{eq:uniform-mu} implies $\lim_{\vartheta_1\to 0} \vartheta_1 \mu(\vartheta_1,\vartheta_2)=0$.
Furthermore, once one has a quantitative estimate on the dependence $(\vartheta_1,\vartheta_2)\mapsto \mu$, 
the assumption $4R \vartheta_1 C_3 <1$ can be formulated as a condition involving $\vartheta_1, \vartheta_2$ and $R$ only.

The proof of Theorem~\ref{thm:qUC_elliptic} is based on ideas developed in \cite{Rojas-MolinaV-13} for the pure Laplacian. The key tool for the proof is a Carleman estimate. For second order elliptic operators there exist plenty of them in the literature, see e.g.\ \cite{KochT-01,RousseauL-12,Rousseau-12}. However, since we are interested in quantitative estimates, the Carleman estimate from Theorem~\ref{thm:Carleman} proved to be useful in this context.
%


\subsection*{Acknowledgment}
I.V. would like to thank Thomas Duyckaerts and 
Matthieu L\'eautaud for discussions concerning the literature on Carleman estimates.

\begin{thebibliography}{17}

\bibitem{BourgainK-05}
J.~Bourgain and C.~E. Kenig.
\newblock On localization in the continuous {A}nderson-{B}ernoulli model in
  higher dimension.
\newblock {\em Invent. Math.}, 161(2):389--426, 2005.


\bibitem{BourgainK-13}
J.~Bourgain and A.~Klein.
\newblock Bounds on the density of states for {Schr\"odinger} operators.
\newblock {\em Invent. Math.},  194(1):{41--72}, {2013}.

\bibitem{BorisovTV}
D.~I. Borisov, M.~Tautenhahn, and I.~Veseli\'c.
\newblock Equidistribution properties of eigenfunctions of divergence form
  operators.
\newblock in preparation.


\bibitem{ButzerSS-88}
P.~L.~Butzer, W.~Splettstoesser, and R.L.~Stens.
\newblock 
The sampling theorem and linear prediction in signal analysis.
\newblock  Jahresber. Deutsch. Math.-Verein. 90:1-70, 1988. 


\bibitem{Carleman-39}
T.~Carleman.
\newblock Sur un probl\`{e}me d'unicit\'{e} pour les syst\`{e}mes
  d'\'equations aux d\'eriv\'ees partielles \`a deux variables ind\'ependantes.
\newblock {\em Ark. Mat. Astr. Fys.}, 26(17):1--9, 1939.

\bibitem{EscauriazaV-03}
L.~Escauriaza and S.~Vessella.
\newblock {\em Optimal Three Cylinder Inequalities for Solutions to Parabolic
  Equations with {L}ipschitz Leading Coefficients}, volume 333 of {\em Contemp.
  Math.}, pages 79--87.
\newblock American Mathematical Society, 2003.

\bibitem{GerminetK-13}
F.~Germinet and A.~Klein.
\newblock A comprehensive proof of localization for continuous {A}nderson models with singular random potentials.
\newblock {\em J. Eur. Math. Soc.}, 15(1):53--143, 2013.

\bibitem{Klein-13}
A.~Klein.
\newblock Unique continuation principle for spectral projections of
  {S}chr{\"o}dinger operators and optimal {W}egner estimates for non-ergodic
  random {S}chr{\"o}dinger operators.
\newblock {\em Commun. Math. Phys.}, 323(3):1229--1246, 2013.

\bibitem{KochT-01}
H.~Koch and D.~Tataru.
\newblock Carleman estimates and unique continuation for second-order elliptic
  equations with nonsmooth coefficients.
\newblock {\em Commun. Pur. Appl. Math.}, 54(3):339--360, 2001.

\bibitem{RousseauL-12}
J.~{Le Rousseau} and G.~Lebeau.
\newblock On {C}arleman estimates for elliptic and parabolic operators.
  {A}pplications to unique continuation and control of parabolic equations.
\newblock {\em ESAIM Contr. Optim. Ca.}, 18(3):712--747, 2012.

\bibitem{Mueller-54}
C.~M{\"u}ller.
\newblock On the behavior of the solutions of the differential equation $\delta
  u = f(x,u)$ in the neighborhood of a point.
\newblock {\em Commun. Pur. Appl. Math.}, 7(3):505--515, 1954.

\bibitem{NakicRT}
I.~Naki\'c, C.~Rose, and M.~Tautenhahn.
\newblock A quantitative Carleman estimate for second order elliptic operators. 
\newblock arXiv:1502.07575 [math.AP], 2015.

\bibitem{NakicTV}
I.~Naki\'c, M.~T\"aufer, M.~Tautenhahn, and I.~Veseli\'c.
\newblock Scale-free uncertainty principles and Wegner estimates for random breather potentials.
\newblock arXiv:1410.5273 [math.AP], 2014.

\bibitem{NakicTTV-prep}
I.~Naki\'c, M.~T\"aufer, M.~Tautenhahn, and I.~Veseli\'c.
\newblock In preparation.

\bibitem{Rojas-MolinaV-13}
C.~Rojas-Molina and I.~Veseli{\'c}.
\newblock Scale-free unique continuation estimates and applications to random
  {S}chr{\"o}dinger operators.
\newblock {\em Commun. Math. Phys.}, 320(1):245--274, 2013.

\bibitem{Rousseau-12}
J.~Le Rousseau.
\newblock Carleman estimates and some applications to control theory.
\newblock In {\em Control of Partial Differential Equations}, Lecture Notes in Mathematics, pages 207--243. 2012.

\bibitem{Vessella-08}
S.~Vessella.
\newblock Quantitative estimates of unique continuation for parabolic
  equations, determination of unknown time-varying boundaries and optimal
  stability estimates.
\newblock {\em Inverse Probl.}, 24(2):023001, 2008.

\bibitem{Ziemer-89}
W.~P. Ziemer.
\newblock {\em Weakly differentiable functions}.
\newblock Springer, New York, 1989.

\end{thebibliography}
\end{document}